\documentclass[12pt,twoside]{amsart}
\usepackage{amsmath}
\usepackage{amssymb}
\usepackage{amscd}
\usepackage{latexsym}
\usepackage{enumerate}
\usepackage[mathscr]{eucal}

\newtheorem{theorem}{Theorem}[section]

\newtheorem{definition}[theorem]{Definition}

\newcommand{\K}{{\mathfrak K}}

\newcommand{\adic}{^{\raise0.05cm\hbox{\scriptsize $\wedge$}}}

\def\arr{\rightarrow}

\newcommand{\downarrowright}[1]{\downarrow
\rlap{\raise0.1cm\hbox{$\scriptstyle{#1}$}}}

\newcommand{\downarrowleft}[1]{\rlap{\kern-0.2cm
\raise0.1cm\hbox{$\scriptstyle{#1}$}}\downarrow}

\newcommand{\uparrowright}[1]{\uparrow
\rlap{\lower0.1cm\hbox{$\scriptstyle{#1}$}}}

\newcommand{\uparrowleft}[1]{\rlap{\kern-0.2cm
\lower0.1cm\hbox{$\scriptstyle{#1}$}}\uparrow}

\title[Radicals and Plotkin's problem concerning geometrically equivalent groups] {Radicals and Plotkin's problem concerning geometrically equivalent groups}

\author[R\"udiger G\"obel]{ R\"udiger G\"obel (Essen)}
\address{Fachbereich 6, Mathematik und Informatik, Universit\"at Essen, 45117 Essen, Germany }
\email{R.Goebel@Uni-Essen.De}
\author[Saharon Shelah]{Saharon Shelah (Jerusalem)}

\address{Department of Mathematics,
Hebrew University, Jerusalem, Israel\\
and Rutgers University, New Brunswick, NJ, U.S.A} \email{e-mail:
Shelah@math.huji.ac.il}

\thanks{The authors are supported by the project No. G
0545-173, 06/97 of the German-Israeli Foundation for Scientific
Research \& Development.
\newline
GbSh 741 in Shelah's list of publications.}

\begin{document}

\begin{abstract}
If $G$ and $X$ are groups and $N$ is a normal subgroup of $X$,
then the $G-$closure of $N$ in $X$ is the normal subgroup
${\overline X}^G = \bigcap \{ \ker \varphi | \varphi :
X\rightarrow G, \mbox{ with } N \subseteq \ker \varphi \}$ of
$X$. In particular, ${\overline 1}^G = R_GX$ is the $G-$radical of
$X$. Plotkin \cite{KM, PPT,P1} calls two groups $G$ and $H$
geometrically equivalent, written $G\sim H$,
 if for any free group $F$ of finite rank
and any normal subgroup $N$ of $F$ the $G-$closure and the
$H-$closure of $N$ in $F$ are the same. Quasiidentities are
formulas of the form $(\bigwedge_{i\le n} w_i = 1 \rightarrow w
=1)$ for any words $w, w_i \ (i\le n)$ in a free group. Generally
geometrically equivalent groups satisfy the same quasiidentiies.
Plotkin showed that  nilpotent groups $G$ and $H$ satisfy
the same quasiidenties  if and
only if $G$ and $H$ are geometrically equivalent. Hence he
conjectured that this might hold for any pair of groups; see the
Kourovka Notebook \cite{KM}. We provide a counterexample.

\end{abstract}

\date{}
\maketitle

In a series of paper, B. I. Plotkin and his collaborators
\cite{PPT,P1,P2,P3} investigated radicals of groups and their
relation to quasiidentities. If $G$ is a group, then the {\it
$G$-radical $R_G X$ of a group $X$} is defined by
$$R_G X = \bigcap \{\ker \varphi ; \varphi: X \arr G \mbox{ any
homomorphism } \}.$$

 Clearly,
$R_G X$ is a characteristic, hence a normal subgroup of $X$. The
radical $R_G$ can also be used to define the  {\it $G$-closure
${\overline U}^G={\overline U}$ of a normal subgroup $U$ of $X$,}
by saying that ${\overline U}/U = R_G (X/U)$. This immediately
leads to Plotkin's definition of geometrically equivalent groups,
see \cite {PPT, P1, P2, P3} and \cite [p. 113] {KM}.

\begin{definition}  Let $G$ and $H$ be two groups. Then $G$ and $H$ are
geometrically equivalent, written $G \sim H$, if for any free
group $F$ of finite rank and any normal subgroup $U$ of $F$ the
$G$- and $H$-closure of $U$ in $F$ are the same, i.e. for any
normal subgroup $U$ we have ${\overline U}^G = {\overline U}^H$.
\end{definition}

It is easy to see that $G \sim H$ if and only if $R_G K = R_H K$
for all finitely generated groups $K$. Plotkin notes that
geometrically equivalent groups satisfy the same quasiidentities.
The well-known notion of quasiidentities relates to
quasivarieties of groups. A {\it quasiidentity } is an expression
of the form
$$w_1 = 1 \wedge \dots \wedge w_n = 1 \arr w = 1 \mbox{ where }
w_i, w \in F \ (i \leq n) \mbox{ are words.}$$ Moreover the
following was shown in \cite{PPT}, see \cite[p.113] {KM}.

\begin{theorem} \label{plot}
\begin{enumerate}
\item[(a)]  If $G \sim H$, and $G$ is torsion-free, then $H$ is
torsion-free.
\item[(b)] If $G, H$ are nilpotent, then $G \sim H$
if and only if $G$ and $H$ satisfy the same quasiidentities.
\end{enumerate}
\end{theorem}

This lead Plotkin to conjecture that two groups might be
geometrically equivalent if and only if they satisfy the same
quasiidentities, see the Kourovka Notebook \cite [p.113, problem
14.71] {KM}. In this note we refute this conjecture.
Clearly there are only countably many finitely presented groups
which we enumerate as the set $\K = \{K_n : n \in w\}$ and let $G
= \prod\limits_{n \in w} K_n$ be the restricted direct product.
Then $G$ satisfies only those quasiidentities satisfied by all
groups and so if $H$ is any group with $G \leq H$,
$G$ satisfies the same quasiidenties as $H$.

R. Camm \cite[p. 68, p. 75 Corollary]{Ca} proved there are
$2^{\aleph_0}$ non-isomorphic, two-generator, simple groups, see
also Lyndon, Schupp \cite[p. 188, Theorem 3.2]{LS}.  So there
exists a 2-generated simple group $L$ which cannot be mapped
nontrivially
 into $G$.
 We consider the pair $G$,  $H = L \times G$ and show the
following:

\begin{theorem} \label{equi}  If $G, H$ and $L$ are as above, $R_G L = L$ and $R_H L = 1$.
 In particular $G$ and $H$ are not geometrically
equivalent.
 Since $G\leq H$ satisfy the same quasiidentities, this is the required counterexample.
\end{theorem}

\begin{proof} Since $L$ is a two-generated
simple group, $L$ is an epimorphic image of a free group of rank
2. So it is enough to prove that $R_G L = L$ and $R_H L = 1$. The
first equality follows since there is no nontrivial homomorphism
of $L$ into $G$. On the other hand, there is a canonical
embedding $L \arr H = L \times G$, so $R_H L = 1$.
\end{proof}

\end{document}